\documentclass{notices}
\usepackage{amsfonts,amssymb,amsmath,amscd}

\usepackage{graphicx}
\usepackage{xcolor}
\usepackage{wrapfig}
\usepackage{CJKutf8}
\usepackage[T1]{fontenc}

\title{
    On the importance of illustration for mathematical research
}

\author{
  R\'emi Coulon
  \affil{
        Univ Rennes, CNRS, IRMAR - UMR 6625, 35000 Rennes, France
    }
  \and
  Gabriel Dorfsman-Hopkins
  \affil{
        Department of Mathematics, Statistics, and Computer Science, St. Lawrence University, Canton, NY
   }
   \and
   Edmund Harriss
   \affil{
        Department of Mathematics and School of Art, University of Arkansas
   }
   \and
   Martin Skrodzki
   \affil{
        Computer Graphics and Visualization, TU Delft, Netherlands
   }
   \and
   Katherine E.~Stange
   \affil{
        Department of Mathematics, University of Colorado Boulder
   }
   \and
   Glen Whitney
   \affil{
        Prison Math Project
   }
}

\begin{document}

\maketitle

\section*{Abstract}

Mathematical understanding is built in many ways. Among these, illustration has been a companion and tool for research for as long as research has taken place. We use the term \emph{illustration} to encompass any way one might bring a mathematical idea into physical form or experience, including hand-made diagrams or models, computer visualization, 3D printing, and virtual reality, among many others.  The very process of illustration itself challenges our mathematical understanding and forces us to answer questions we may not have posed otherwise.  It can even make mathematics an experimental science, in which immersive exploration of data and representations drive the cycle of problem, conjecture, and proof.   Today, modern technology for the first time places the production of highly complicated models within the reach of many individual mathematicians. Here, we sketch the rich history of illustration, highlight important recent examples of its contribution to research, and examine how it can be viewed as a discipline in its own right.

\section*{Introduction} 

\begin{figure}
    \centering
    \includegraphics[width=0.95\linewidth]{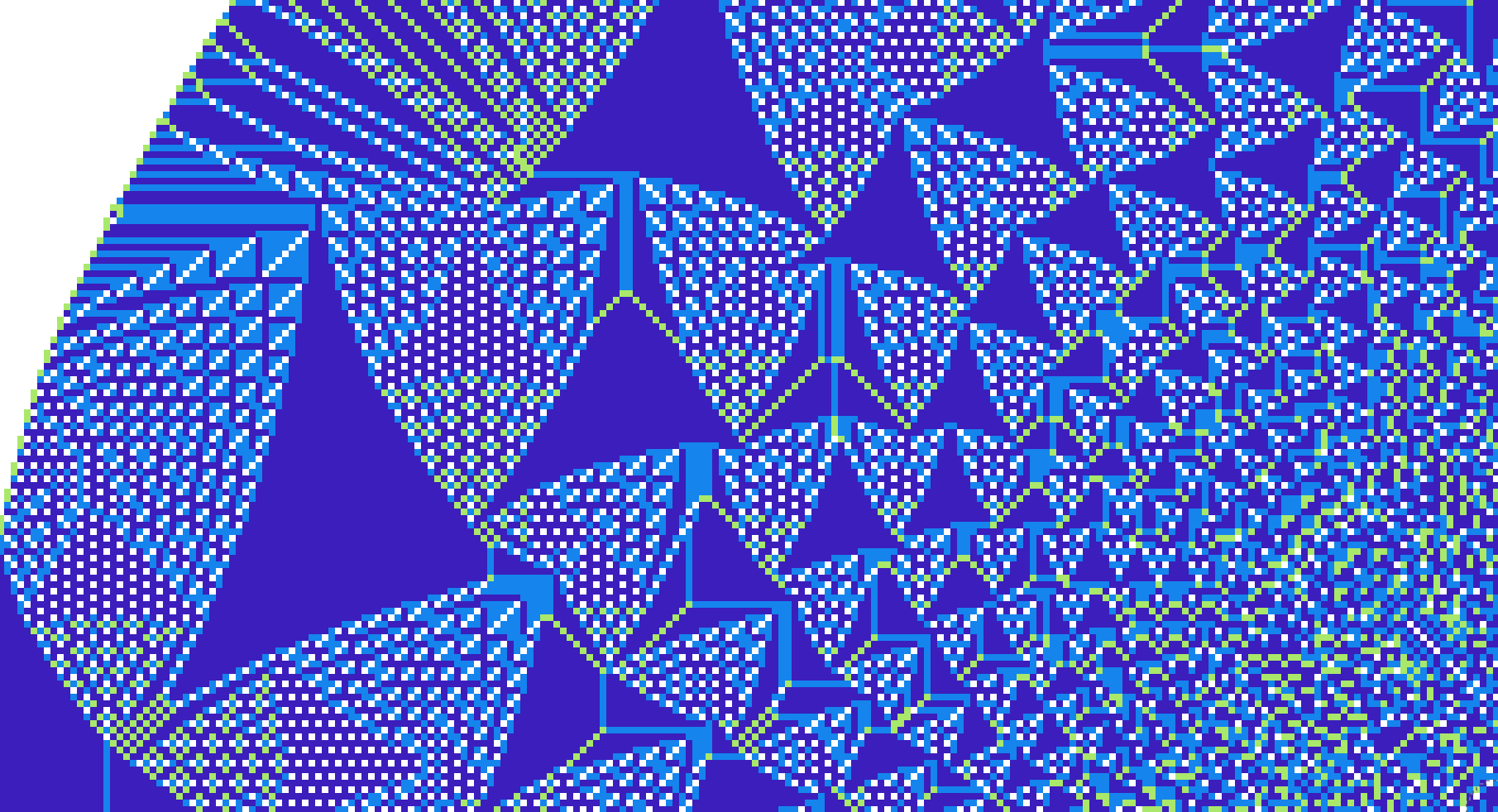}
    \caption{Detail of a $\mathbb{Z}^2$-lattice sandpile stabilized from a starting configuration with 500 million grains at the origin.  Whenever any position has at least 4 grains, it distributes one to each neighbour.  Nodes in the stable configuration have 0, 1, 2 or 3 grains (indicated by colour).  Image by Stange.}
    \label{fig:sandpile}
\end{figure}

In the last decade, it has become increasingly possible for researchers to augment their experience of abstract mathematics with sensory exploration: 3D-printed or CNC-milled models, the ability to walk through seemingly impossible physical spaces with virtual reality, and the potential to explore high-dimensional mathematical spaces through computer visualisation all provide such opportunities, among other media. Now much more than simply an aid to understanding, these tools have reached a level of sophistication that makes them indispensable to many frontiers of mathematical research. To preview one particular case recounted below, the tantalizing structure visible in Figure \ref{fig:sandpile} (and many others like it) led to conjectures and proofs that would likely otherwise have been inaccessible. 

The list of examples of research driven by illustration is rapidly expanding in recent years. 
We use the term \emph{illustration} to encompass any way one might bring a mathematical idea into physical form or experience, including hand-made diagrams or models, computer visualization, 3D printing, or virtual reality, among many others. 
We will discuss instances of this interplay including breakthroughs in such areas as representation theory (such as Knutson and Tao's honeycombs~\cites{honeycombAMS,honeycomb1,honeycomb2}), self-organized criticality (work of Levine, Pegden, Smart and others on Abelian sandpiles~\cite{LPS17} as illustrated in Figure \ref{fig:sandpile}), geometry (Coiculescu and Schwartz on Thurston's \emph{Sol} geometry~\cite{coiculescu2022sol}), and most recently the solution of the einstein problem with the hat monotile~\cite{smith2023aperiodic} and its chiral version, the spectre~\cite{smith2023chiral}. 
These solutions came after many years of physical exploration by Dave Smith, which revealed especially interesting properties of the hat tile.
Note that these papers have already sparked several follow-up works~\cites{reitebuch2023direct,baake2023dynamics,socolar2023quasicrystalline}).

Illustration is beginning to find a home at programs like the special semester in \emph{Illustrating Mathematics} in Fall 2019 at the Institute for Computational and Experimental Research in Mathematics (ICERM) and the Institute for Advanced Study (IAS)/Park City Math Institute (PCMI) virtual program in Summer 2021,\footnote{See \url{http://illustratingmath.org/} for links to these two programs, along with many other resources.} and a community is forming around many modern tools.
\begin{CJK*}{UTF8}{bsmi}
Of course, the importance of illustration to research is not new:  abstraction was linked to plane diagrams in the work of the ancients, including Euclid's \emph{Elements} or the Chinese treatise \emph{The Nine Chapters on the Mathematical Art} (九章算術). Precise three-dimensional models were produced by skilled artisans in the 19th century, notable examples of which remain in the collections at the Institute Henri Poincaré\footnote{\url{https://patrimoine.ihp.fr/}}~\cite{Villani:OM}
\end{CJK*}
or G\"ottingen University\footnote{\url{https://sammlungen.uni-goettingen.de/sammlung/slg_1017/}}, among many other institutions. When computer visualization first became widely available in the 1980's, the Geometry Center was founded at the University of Minnesota, with a mission to exploit these new tools on behalf of mathematics. 
But we are now at another cusp:  modern technological tools have recently made 3D models and virtual reality widely available, and computation and computer visualization is more accessible and more powerful than ever.  We can now collect huge mathematical datasets and examples, and it has become urgent to develop intuitive ways to interact with this data.  

Making full use of modern tools is not without its challenges:  beyond the obvious technical challenges and software learning curves, there are important questions about how an illustration, much like a statistic or an experiment, can subtly mislead the researcher, or miss the essential mathematical pattern sought. Researchers often individually reinvent the necessary skill sets as they seek to advance their own projects, and these projects are pushing the boundaries of the possible. 
But by building a discipline around this enterprise, we can develop its full potential to advance mathematical research.

\section*{Some highlights from the history of mathematical illustration}

\begin{figure}
    \includegraphics[width=1.\linewidth]{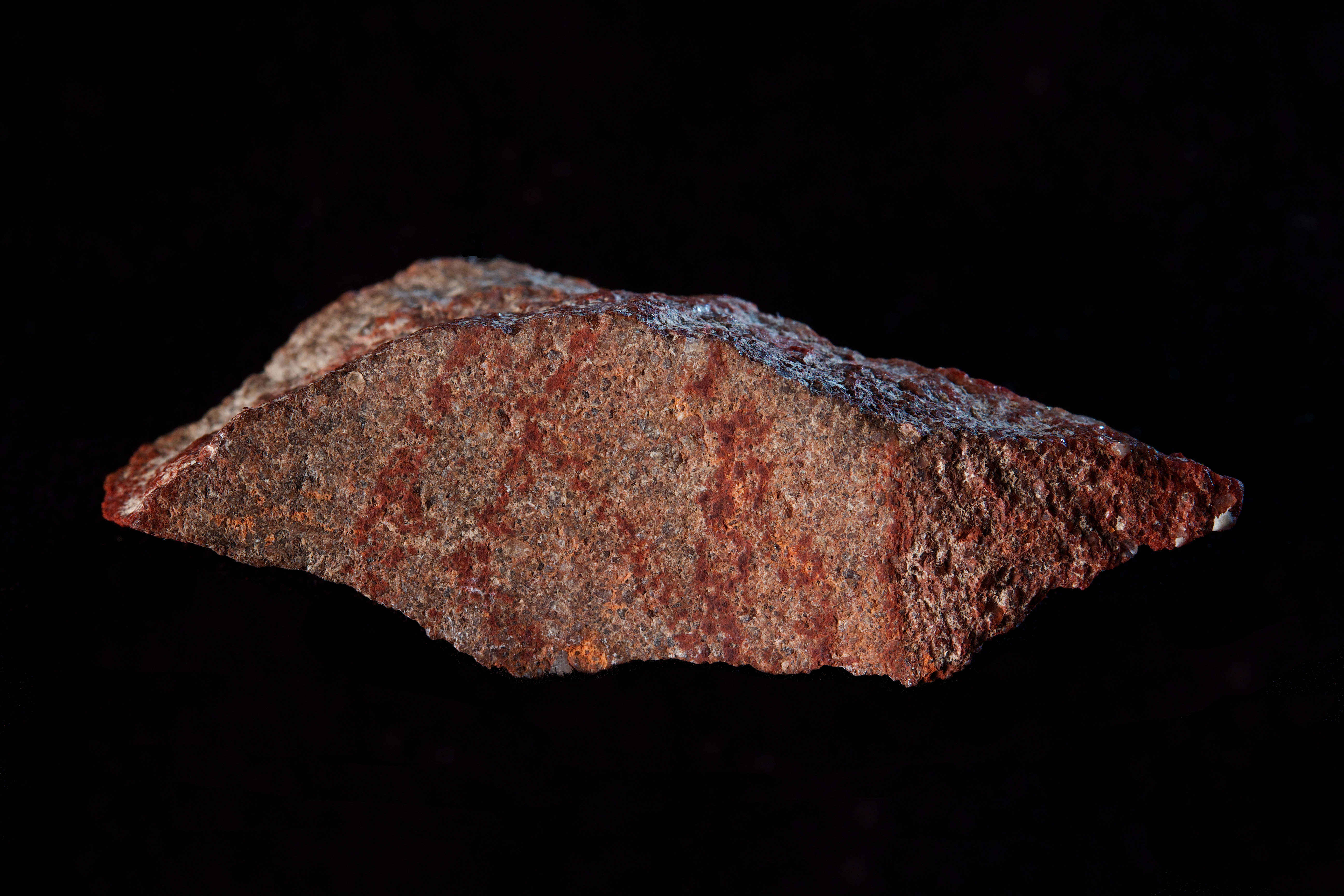}
    \caption{An abstract ochre drawing on silica from the 73,000 year old layers of the Blombos Cave, studied in~\cite{Henshilwood:ADF7LBCSA}. 
    Image by Craig Foster courtesy of Professor Christopher Henshilwood.}
    \label{fig:BlombasStone}
\end{figure}

\begin{figure}
    \includegraphics[width=1.\linewidth]{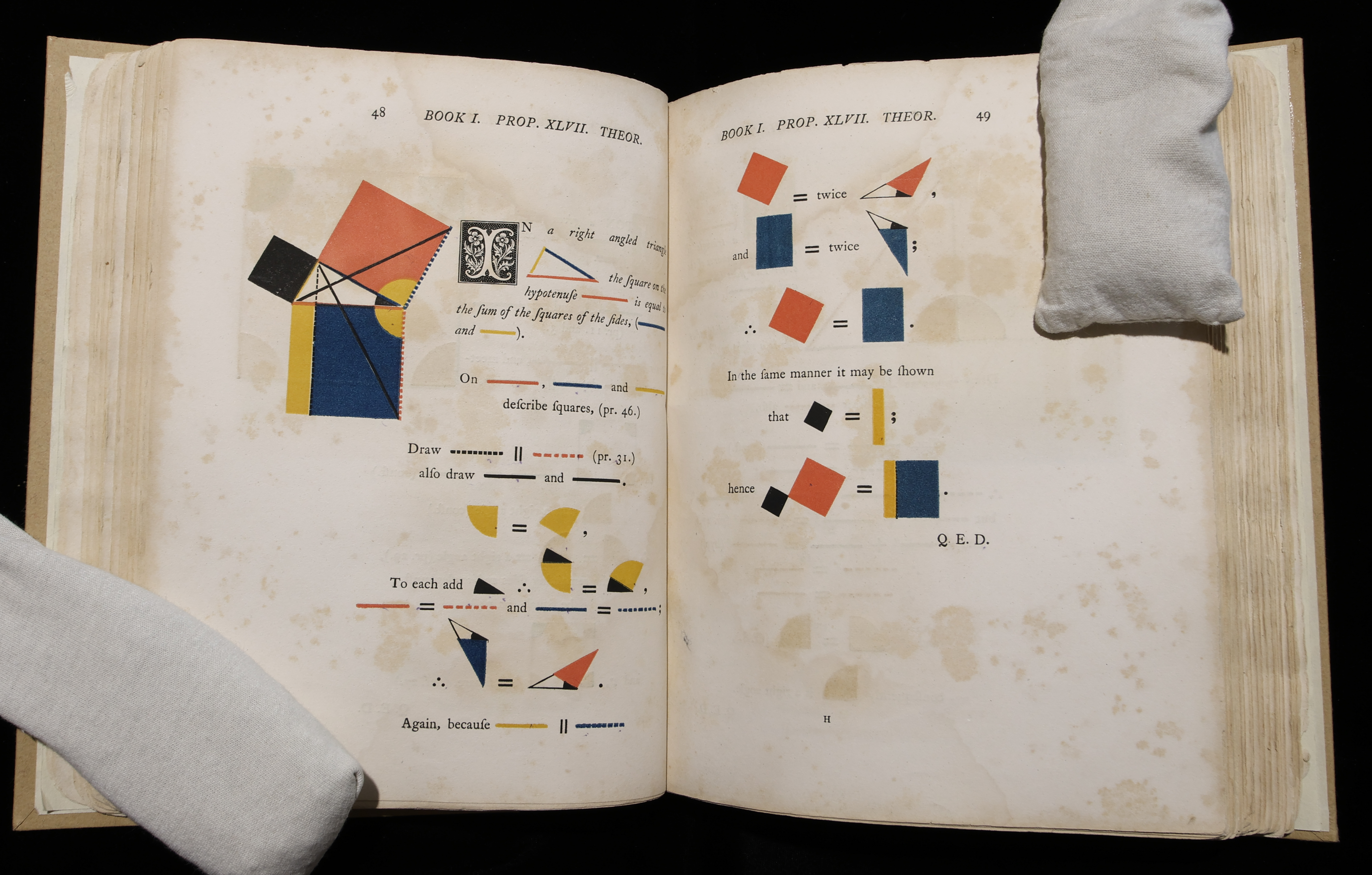}
    \caption{The proof of Book 1 Proposition 47, half of Pythagoras' theorem, in an Oliver Byrne's edition from 1847. The figure shown here had appeared in printed copies of the \emph{Elements} since the fifteenth century, but Byrne's rendition links it tightly and visually to the proof alongside. Image Credit: Harriss and University of Arkansas Library.}
    \label{fig:Pythagoras}
\end{figure}

Illustration of mathematics goes back as far as mathematical ideas themselves. In fact, some of the earliest evidence we have for abstract thinking comes from human-made designs, for example the cross-hatched carvings in Blombos Caves in South Africa, potentially from 73,000 years ago~\cite{Henshilwood:ADF7LBCSA}. 
A little more recently, the middle-eastern tradition of geometry presented in Euclid's \emph{Elements}~\cite{Euclid:EE}
provides a structural link between statements deduced from axioms, and figures made with straight edge and compass.  These two tools provide a physical realization of the two key objects (straight lines and circles) described by the axioms. Euclid's diagrams give a map to help follow (or discover!) the chain of deduction in a proof.  Conversely, the proof validates the image (which could otherwise mislead by error or the selection of a non-generic example).
This approach leads at the conclusion of Book 1 to a proof of the Pythagorean theorem; see Figure~\ref{fig:Pythagoras}.

\begin{figure}
    \includegraphics[width=1.\linewidth]{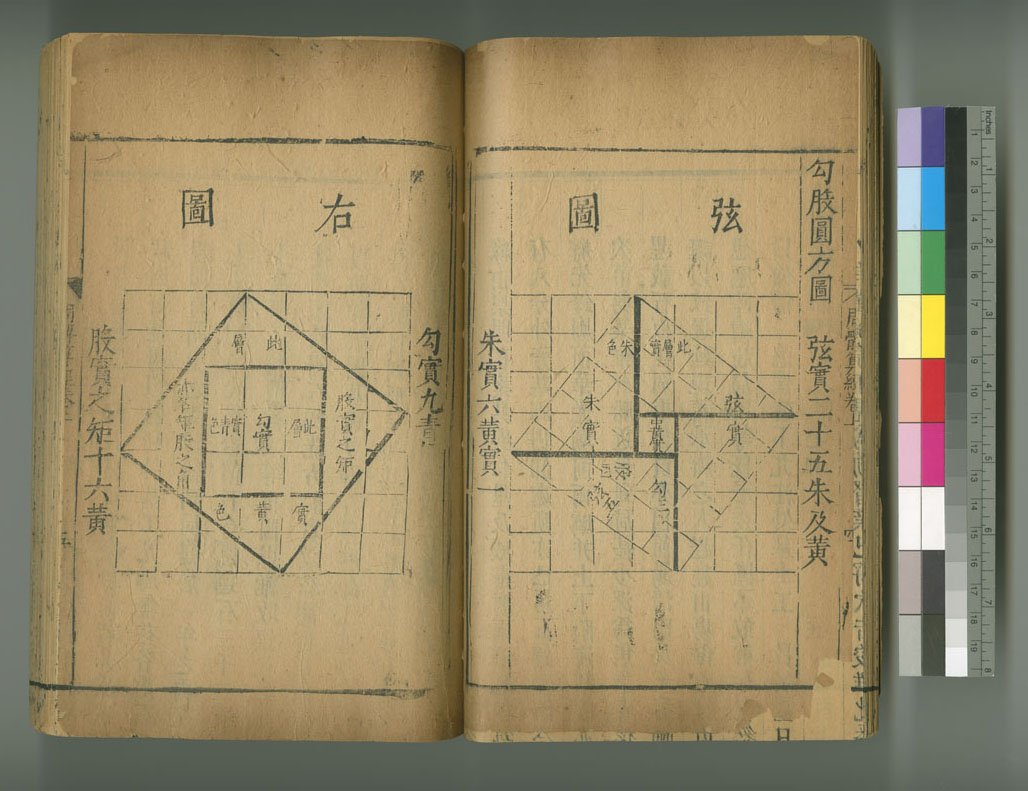}
    \caption{Two pages from the \textit{Arithmetical Classic on the Gnomon and the Circular Paths of Heaven} (周髀算經), a Chinese work on astronomy and mathematics showing a proof of the right triangle 勾股 \ (Gougu) theorem. Image credit: Zhou bi suan jing : shang xia juan, volume 1 / Zhao, Junqing; Zhen, Luan ; Li, Chunfeng et al.; Hu, Zhenting; Bao, Shan = 周髀算經 : 上卷 / 趙君卿注,甄鸞重述,李淳風等注釋,胡震亭提辭,鮑山較, circa 1621-1627. Smith Chinese C-13, folio 4-5. Smith-Plimpton East Asian Collection. Rare Book and Manuscript Library, Columbia University Library.}
    \label{fig:ChineseTriangles}
\end{figure}

\begin{CJK*}{UTF8}{bsmi}
In Chinese mathematics, this theorem is the 勾股 \ (Gougu) theorem. In the classic \emph{Nine Chapters on the Mathematical Art} (九章算術)~\cite{Shen:NCMCC}
, it plays a key role in applying the arithmetical mathematics of the text to geometric problems, for example in 
measuring altitude. The Chinese tradition also gives an elegant visual proof of the result by rearranging triangles, as in Figure \ref{fig:ChineseTriangles}. 
\end{CJK*}

\begin{figure}
    \includegraphics[width=1.\linewidth]{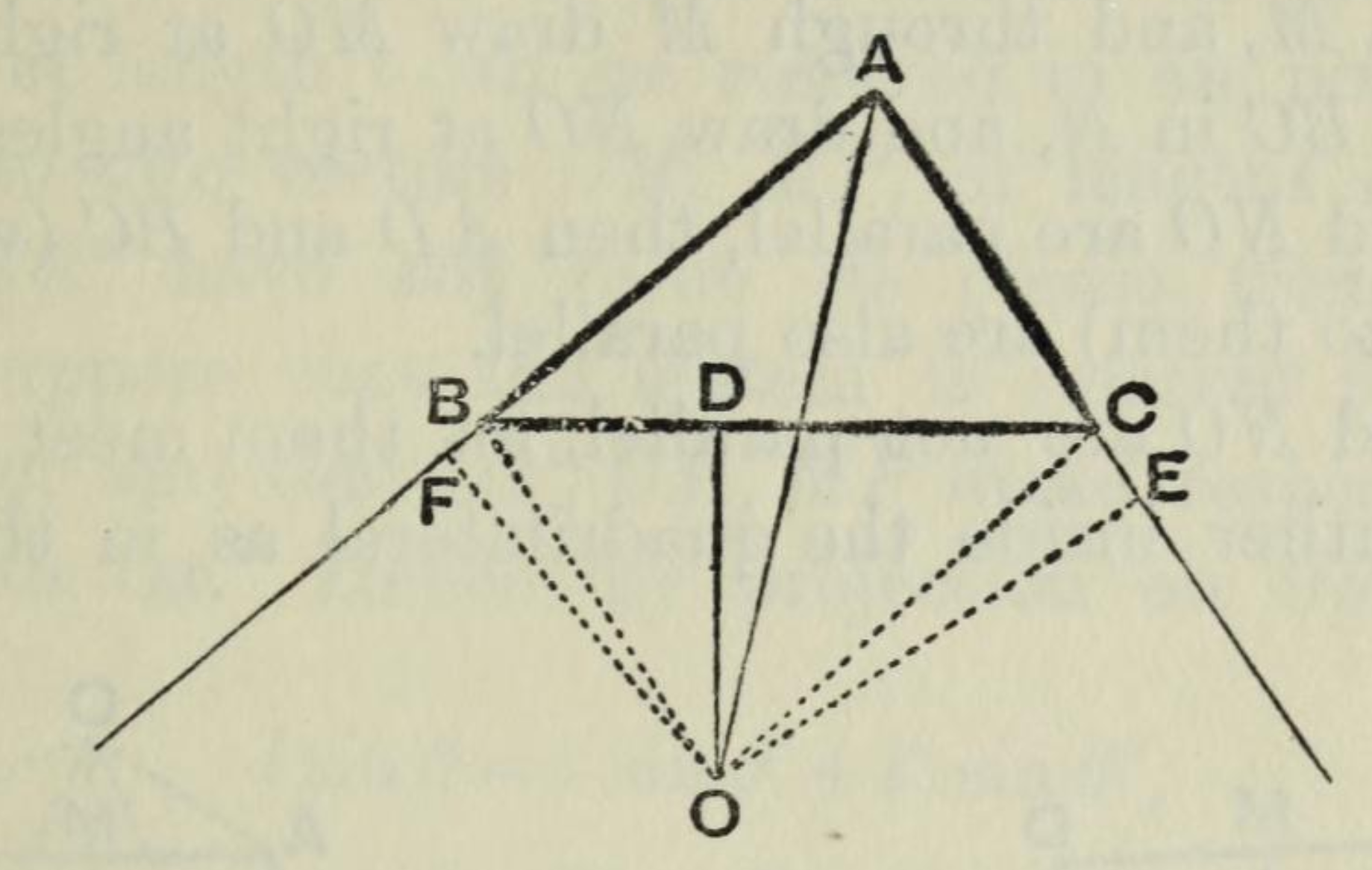}
    \caption{This diagram from Ball's 1882 \emph{Mathematical Recreations and Essays} subtly misleads the reader. In reality, either $E$ lies between $A$ and $C$ or $F$ lies between $A$ and $B$. Image credit: University of Illinois Library/Internet Archive.}
    \label{fig:BallFallacy}
\end{figure}

Although the Chinese proof is not considered rigorous by modern standards, Euclid was also criticized by 
Bertrand Russell when he wrote ``A valid proof retains its demonstrative force when no figure is drawn, but very many of Euclid’s earlier proofs fail before this test.''~\cite{Russell:TE}. This criticism reveals one of the challenges of mathematical illustration.
 A powerful example comes from a well-known ``proof,'' apparently originally due to W. W. Rouse Ball, that all triangles are equilateral.\footnote{For a modern presentation of this fallacy and a lively discussion of its flaws, see Joel Hamkins' essay at \url{https://jdh.hamkins.org/all-triangles-are-isosceles/}.} In this case a subtly misleading diagram leads to a false conclusion following an otherwise entirely correct argument. Briefly, starting with no constraints on the triangle $ABC$, Ball produces points $O$, $E$, and $F$ by a construction depicted in Figure \ref{fig:BallFallacy}. He establishes (without error) that $AF=AE$ and $FB=EC$. The diagram then makes it ``clear'' that $AB = AF-FB$ must equal $AC = AE-EC$. Hence the triangle is isosceles, and since the labels of $A$, $B$, and $C$ were arbitrary, it ``must'' in fact be equilateral!
 
 Disallowing these particular subtle errors requires axioms capturing the meaning of ``between,'' that took considerable work by David Hilbert to formulate~\cites{Hilbert:FG, Lee:AG}. 
 Armed with these axioms, one can show that, except when the triangle \emph{is} isosceles, Ball's construction perfectly executed will produce either $E$ within segment $AC$ or $F$ in segment $AB$ (but not both), destroying the  final subtraction step.

\begin{figure}
    \includegraphics[width=1.\linewidth]{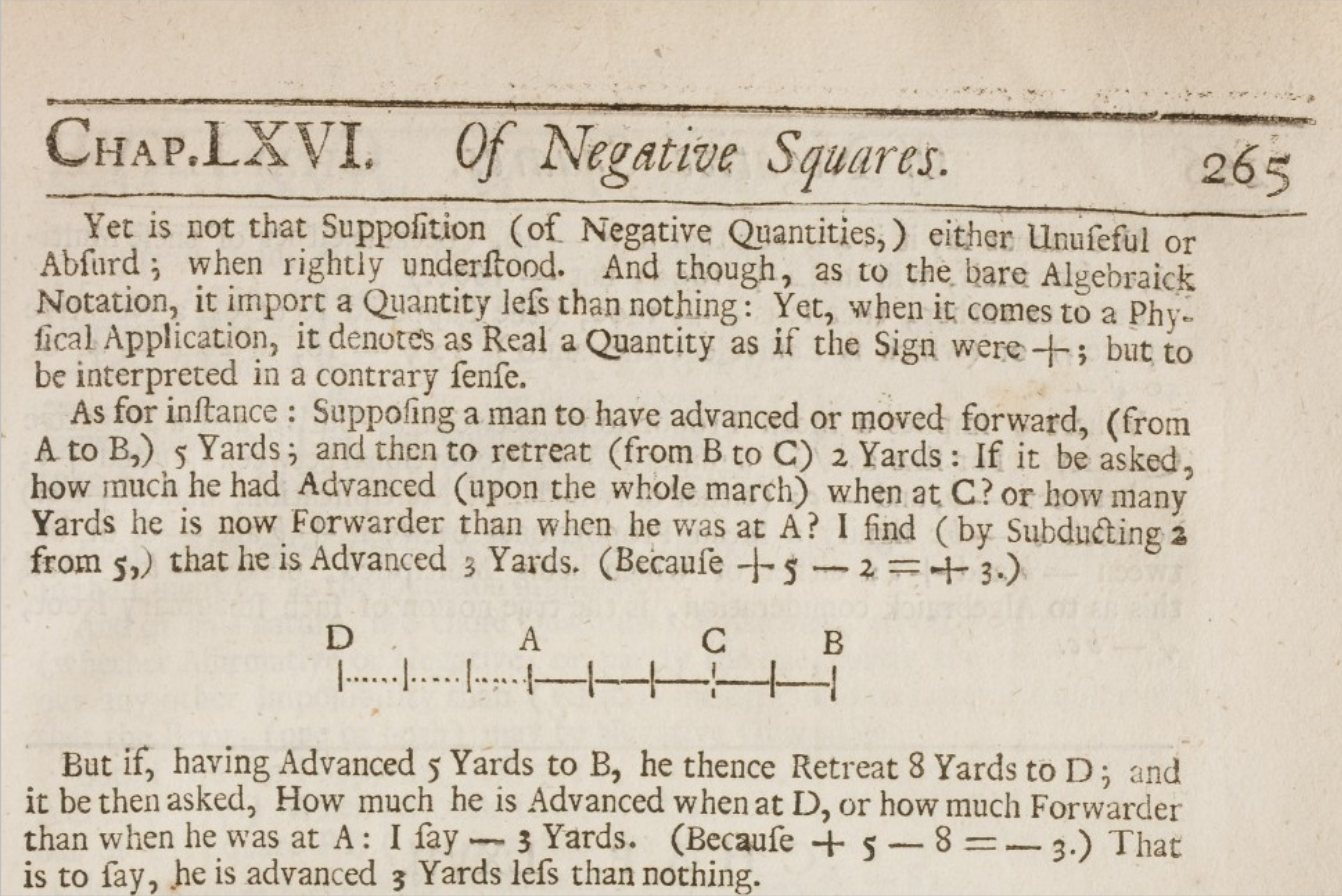}
    \caption{An early depiction of the now-familiar number line, from Wallis' 1685 \emph{A Treatise of Algebra}; image credit: Max Planck Institute for the History of Science, Library.}
    \label{fig:wallisLine}
\end{figure}

A related pitfall -- when a good illustration, overused, can become a pair of blinders -- is illustrated by the following example.
In the \emph{Elements}, the concept of number is based on the concept of length. So the squares in the Pythagorean theorem are actual squares (the area of which are equal), not squared numbers. In the 11th century algebra treatise of Omar Khayyam, although he gives solutions to equations with higher powers than three, he also states: ``Square-square, which, to the algebraists, is the product of the square by itself, has no meaning in continuous objects. This is because how can one multiply a square, which is a surface, by itself? Since the square is a two-dimensional object \ldots\@ and two-dimensional by two-dimensional is a four dimensional object. But solids cannot have more than three dimensions.''~\cite{Khayyam:AOK}. The relation between number and length was also an important factor in the European reluctance to consider negative numbers. A line, after all, cannot have negative length. In contrast, negative quantities are used freely in the \emph{Nine Chapters}, where arithmetic is the foundational idea, with geometry built from it. In Europe the development of the number line, starting with John Wallis, gave an alternative illustration of number (see Figure \ref{fig:wallisLine}) with the capacity to include negative quantities as numbers in their own right~\cite{Schubring:CBGRI}.

Powers and negative numbers are but two examples of a productive pattern of mathematics developing from the tension between illustration and symbolic idea. The study of complex numbers advanced significantly with the concept of the complex plane, and then allowed a new algebraic approach to the geometry of the plane. Both quaternions and matrices were developed to try to extend that understanding to higher dimensions~\cites{Knott:HVM,Grattan-Guiness:MT}. 
In the case of real numbers, although the symbolic ideas would refine the illustrations needed, it was not until the late nineteenth century when fully symbolic definitions were developed, such as Dedekind cuts and Cauchy sequences. At that point the need for illustrations as foundational objects was removed, although the potential for developing intuition and challenging what might be done with the concepts remained~\cite{Fowler:DT}. 

Projective geometry, first developed (as perspective) by artists as a tool to create realistic images, provided one such challenge. These ideas were explored mathematically by Johannes Kepler, G{\'e}rard Desargues and Blaise Pascal. In the early nineteenth century, perspective was developed by Gaspard Monge into ``descriptive geometry'' for the training of engineers in constructing forts and later developed and axiomatized in the foundational work by Jean-Victor Poncelet~\cites{Andersen:GA,Kemp:SOTWFBS}. 
In turn this work would be key in establishing models for non-euclidean geometry, explored axiomatically by Nikolai I.\@ Lobachevsky and J{\'a}nos Bolyai~\cite{Shenitzer:HNGECGS}. 
In this case it was such models, themselves illustrations, that convinced mathematicians of the existence and interest of the non-euclidean geometries. 

Projective geometry also spurred the study of algebraic geometry. In the late nineteenth century an industry 
emerged to reveal the surfaces constructed in this field and their properties, such as cone singularities and embedded straight lines. One pioneer 
was Alexander Brill, a student of Alfred Clebsch with a degree in architecture. Following the work of Peter Henrici (another student of Clebsch), Brill made sliceform paper models of surfaces. He later worked with Felix Klein in Munich to set up a laboratory for the design and production of mathematical 
objects.  
This lab grew into a company that, when it was taken over by Martin Schilling in 1911, had a catalogue of over 400 models. His work combined deep mathematical understanding with a knowledge of printing and construction from his family business~\cites{Polo-Blanco:THGM, OConnor:AWB, Friedman:MMF12C}. 

\begin{figure}
    \includegraphics[width=1.\linewidth]{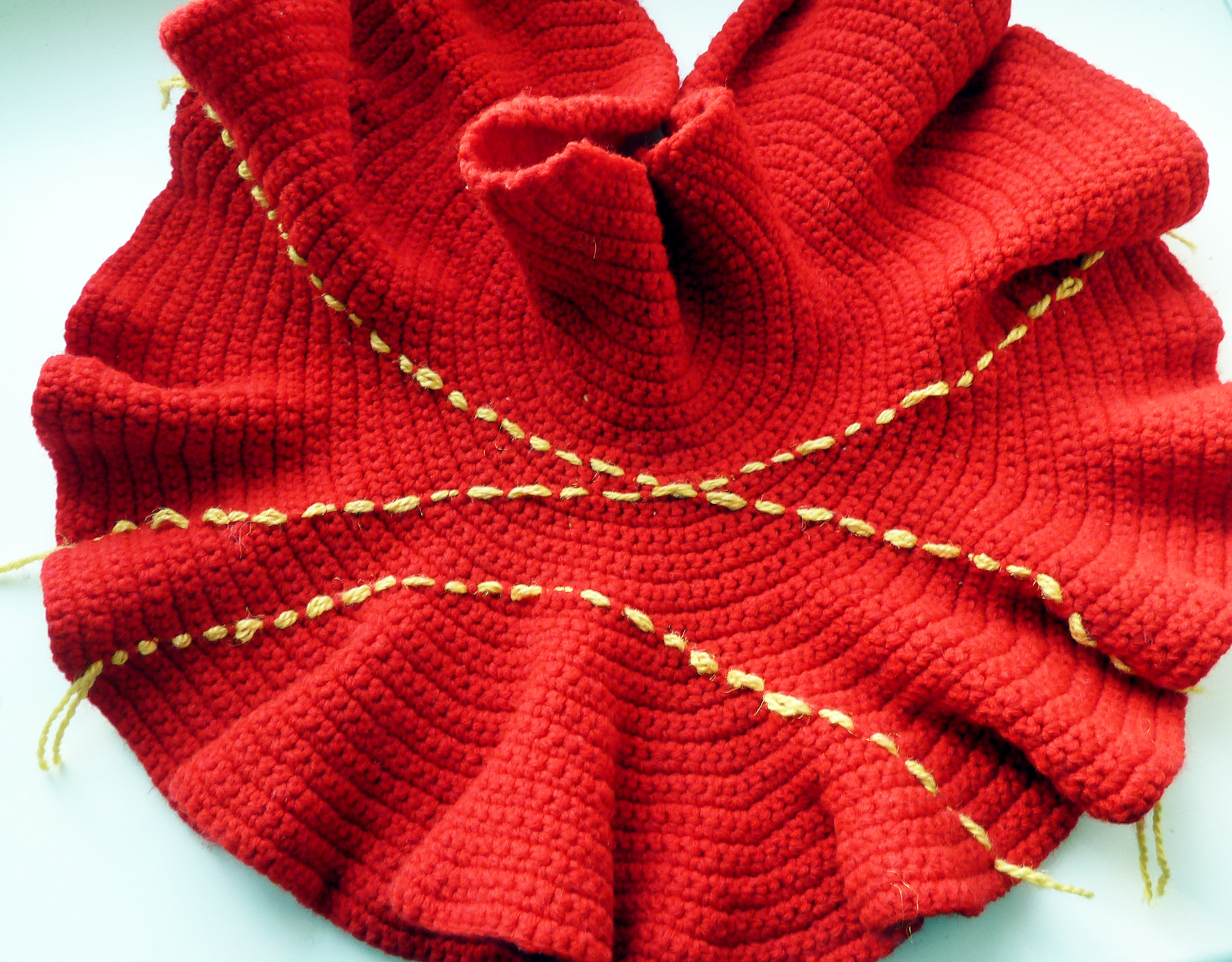}
    \caption{There can be two parallels through a point not on a line, in the (crocheted) hyperbolic plane. Photo credit: Taimiņa.}
    \label{fig:hypercrochet}
\end{figure}

The need to combine mathematical knowledge with fabrication techniques is also highlighted by a story of missed opportunity: how to make physical patches of hyperbolic planes. In addition to his disk model (often called the Poincaré disk model), Eugenio Beltrami also attached together strips of paper to approximate the surface. Other examples used paper polygons connected to make a sort of hyperbolic ``soccer ball.'' These paper models are often fragile, and the rigidity of the paper means that it cannot change its local geometry; thus such models are crude approximations. Roughly a century later, Daina Taimiņa realised that crocheting could produce far more resilient surfaces, with local stretching that meant the negative curvature was more smoothly distributed~\cite{Taimina:CAWHP}. An example of this medium of representation is shown in Figure \ref{fig:hypercrochet}. In fact, similar techniques had been used to create ruffles in scarves and skirts for decades. If the methods of fiber arts had earlier been considered seriously and not dismissed as ``work for women,'' researchers could have had the opportunity to handle robust hyperbolic planes far sooner. 

\section*{The incredible potential for mathematical illustration}

\begin{figure}
    \includegraphics[width=1.\linewidth]{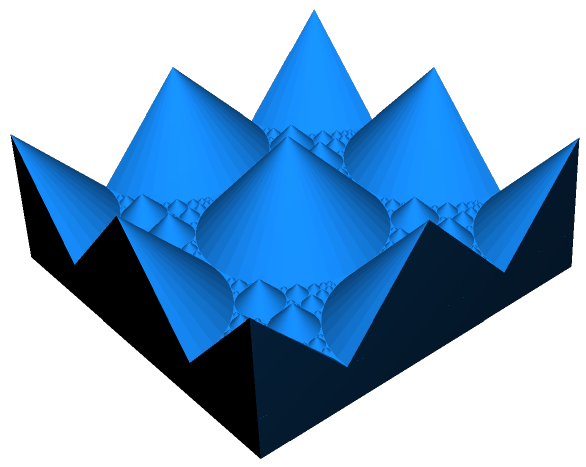}
    \caption{Set of integer superharmonic matrices in the space of all symmetric real matrices. Image by Stange.}
    \label{fig:gamma}
\end{figure}

Turning to recent developments, the work of Lionel Levine, Wesley Pegden, and Charles K.\@ Smart provides an excellent example of the value of illustration as a research tool. Their \emph{Annals of Mathematics} paper \emph{The Apollonian structure of integer superharmonic matrices}~\cite{LPS17} was motivated by the study of Abelian sandpiles on $\mathbb{Z}^2$:  place a large number $N$ of sand grains at the origin, and allow any position with at least four grains to distribute those grains equally to its four neighbours.  The stable configuration that results from this simple system displays impressive large-scale structure that can be discovered through computer visualization (see Figure~\ref{fig:sandpile}). Especially striking is the vivid visual impression that the structure continues to refine at larger $N$ toward a continuum scaling limit, which was proven earlier by Pegden and Smart~\cite{PS13}. To describe the PDE governing this process, the individual periodic tilings in the regions of the limit must be understood.  They are each governed by an \emph{integer superharmonic matrix}.  Levine, Pegden, and Smart generated a picture of the set of integer superharmonic matrices, and were astonished to see the familiar fractal structure of an Apollonian circle packing (Figure~\ref{fig:gamma}).
Each circle of the packing was associated to a periodic pattern appearing in the scaling limit. Through extensive computer investigation, the authors were able to determine the intricate recursive relationships between the patterns for circles generated from one another (`ancestors' overlap and merge to form `descendent' patterns according to complicated rules).  These recursions led to a difficult inductive proof that the set did indeed have the Apollonian structure evident in experiments.  The development of these results provide a perfect example of the role illustration can play in the cycle of conjecture, theorem, and proof.  Without the data available through large-scale computer experimentation and the ability to explore it visually, the question of the scaling limit may not have been raised at all, and the recursive proof of their main result would likely not have been discovered.

\begin{figure*}[ht]
    \centering
    \includegraphics[width=0.9\textwidth]{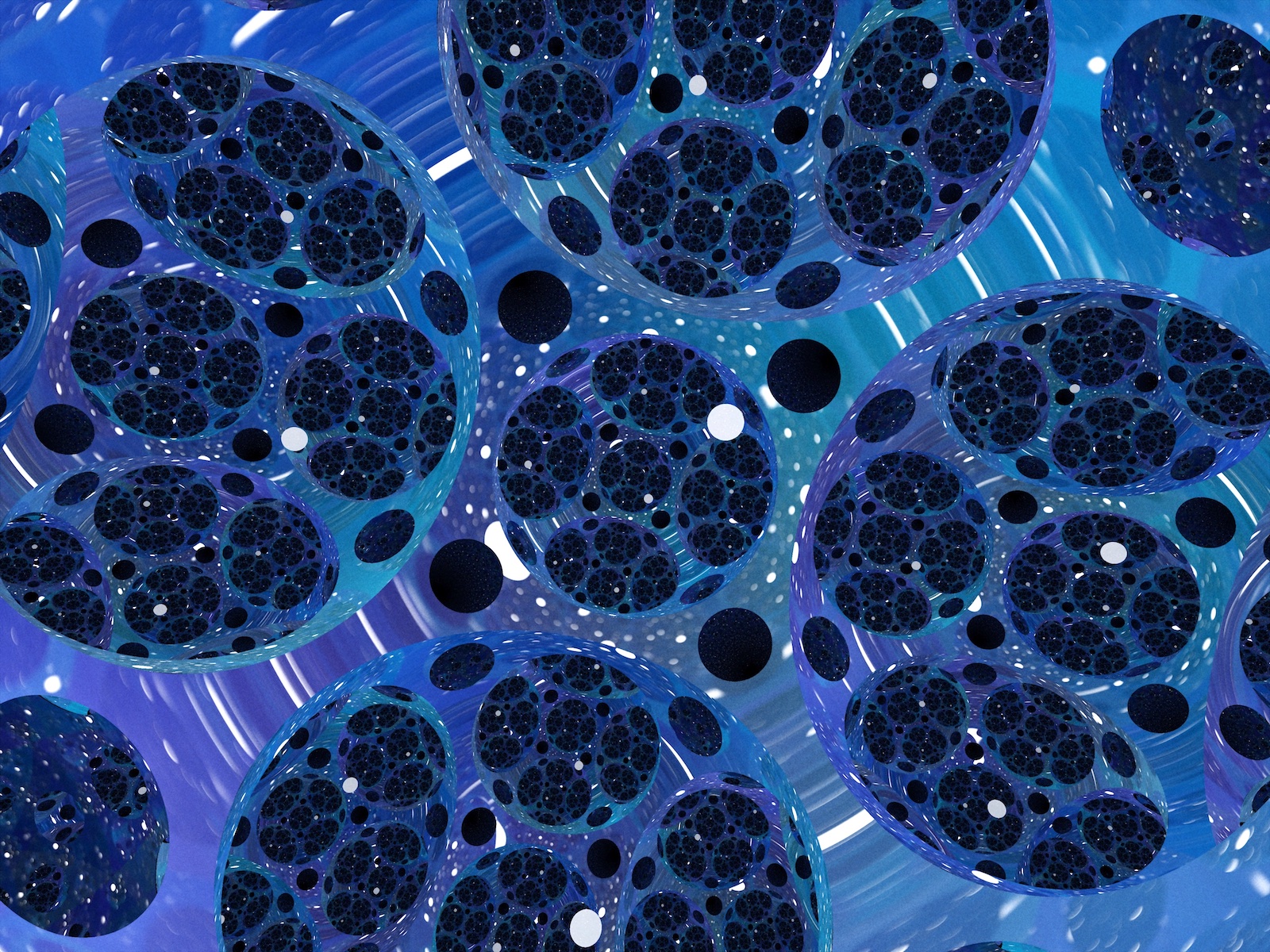}
    \caption{In-space view of a finite-volume hyperbolic $3$-manifold lit by a single white light.
    Image by Coulon, Matsumoto, Segerman and Trettel~\cite{Coulon:2020ac}.}
    \label{fig:hyp}
\end{figure*}

Another area where research is intertwined with illustration is in the study of William Thurston's geometrization conjecture~\cite{geoConjecture}, proved by Grigori Perelman~\cites{perelman1,perelman2,perelman3}.
This key tool in our understanding of $3$-manifolds implies, for instance, the famous Poincaré conjecture. 
Geometrization states that any compact \emph{topological} $3$-manifold can be cut into finitely many pieces, each of which carries a \emph{geometric} structure. 
There are eight possible such structures, known as \emph{Thurston geometries}.
Some of them are rather familiar to mathematicians, such as the 3-dimensional Euclidean and hyperbolic spaces or the $3$-sphere.

Despite the fact that Thurston's geometries have been intensively studied, the more exotic geometries such as \emph{Nil} and \emph{Sol} still defy our ``Euclidean-grown'' spatial intuition.  Keeping in mind the well-established power of our physical and visual intuition to aid geometrical research, R{\'e}mi Coulon, Elisabetta Matsumoto, Henry Segerman, and Steve Trettel developed virtual reality software to immerse the user in any of the eight Thurston geometries~\cite{Coulon:2020ac} (see Figure~\ref{fig:hyp}).  
Besides building much-needed intuition for these spaces, the development of the software itself raised mathematical questions.  
To build their virtually rendered Thurston geometries, these researchers use raymarching techniques which require computation of distances between objects.  But, for example, there is no closed formula for the distance in \emph{Nil} or \emph{Sol}!  Thus, the development of the algorithms themselves becomes a mathematical result in its own right.

\begin{figure*}
    \includegraphics[width=\textwidth]{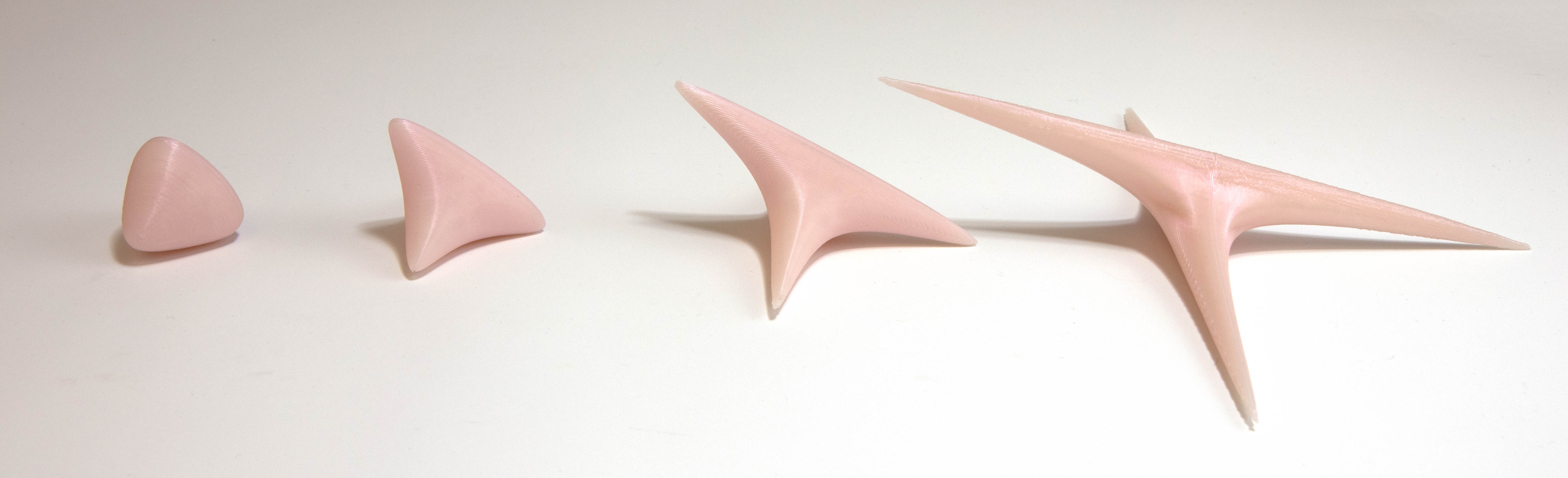}
    \caption{3D printed models of the spheres in Sol produced during the ICERM program ``Illustrating Mathematics'' in Fall 2019.
    Models by Coulon, Image Credit: Harriss.}
    \label{fig:sol spheres}
\end{figure*}

Work on Thurston's geometries has very often been closely tied with illustration. For example, the study of \emph{Spheres in Sol} by Matei P.\@ Coiculescu and Rich Schwartz in \emph{Geometry and Topology} (positively) answers an old open question, whether metric spheres in \emph{Sol} are homeomorphic to $S^2$~\cite{coiculescu2022sol}. Each step of the proof was found after numerous graphical experiments, and 3D printing brings yet another perspective (see Figure~\ref{fig:sol spheres}).

\begin{figure}
    \includegraphics[width=1.\linewidth]{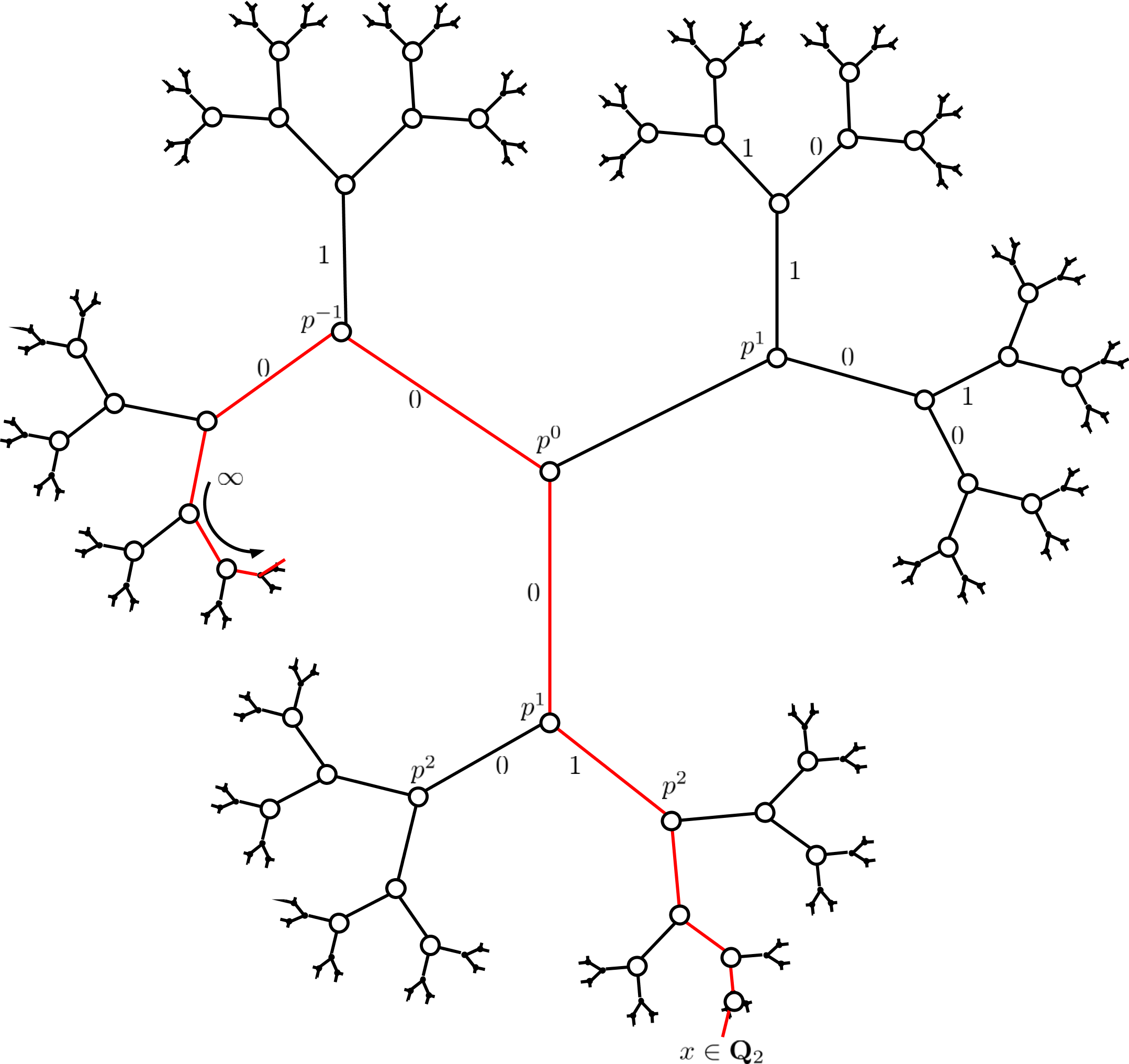}
    \caption{The Bruhat-Tits tree of $\mathbb{Q}_2$, with geodesics~\cite{HMMS17}*{Figure~1}.}
    \label{fig:bruhatTits}
\end{figure}

\begin{figure}
    \includegraphics[width=1.\linewidth]{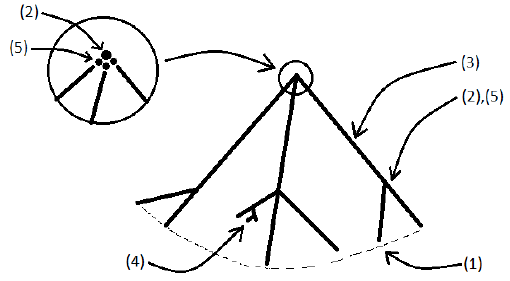}
    \caption{An early illustration of the adic unit ball by Scholze~\cite{sch12}*{Example~2.20}.}
    \label{fig:adicBall}
\end{figure}

 For an example at the intersection of algebraic geometry and number theory, a few key illustrations have helped drive developments in the field of $p$-adic analytic geometry.  At the same time, illustrating the $p$-adic analogs of complex analytic manifolds presents unique challenges, not the least of which is the fact that the $p$-adic numbers themselves are topologically a Cantor set.   Nevertheless, clever and meaningful illustrations of $p$-adic analogs to the complex upper half-plane and complex unit disk have proved incredibly fruitful.
An illustration of Vladimir Drinfeld's \textit{$p$-adic upper half plane} as tubular neighborhoods of Bruhat-Tits trees (Figure~\ref{fig:bruhatTits}) clarified the behavior of the action of $\operatorname{GL}_2(\mathbb{Q}_p)$ by M\"obius transformations.  Understanding this action was instrumental in the work of David Mumford~\cite{mum72}, Drinfeld~\cite{drin76}, in the construction of $p$-adic analytic uniformization of elliptic curves (reflecting the famous complex analytic uniformization of elliptic 
curves as quotients of the complex upper half plane).  Similarly, Peter Scholze's illustrations of the \textit{adic unit ball} (Figure~\ref{fig:adicBall}) provide access to the foundational geometric construction in his theory of perfectoid spaces~\cite{sch12}.  The act of illustrating the central geometric objects of $p$-adic analysis has proven both beneficial and uniquely challenging, demanding a systematic and critical approach.

\begin{figure}
    \centering
    \includegraphics[width=1.0\linewidth]{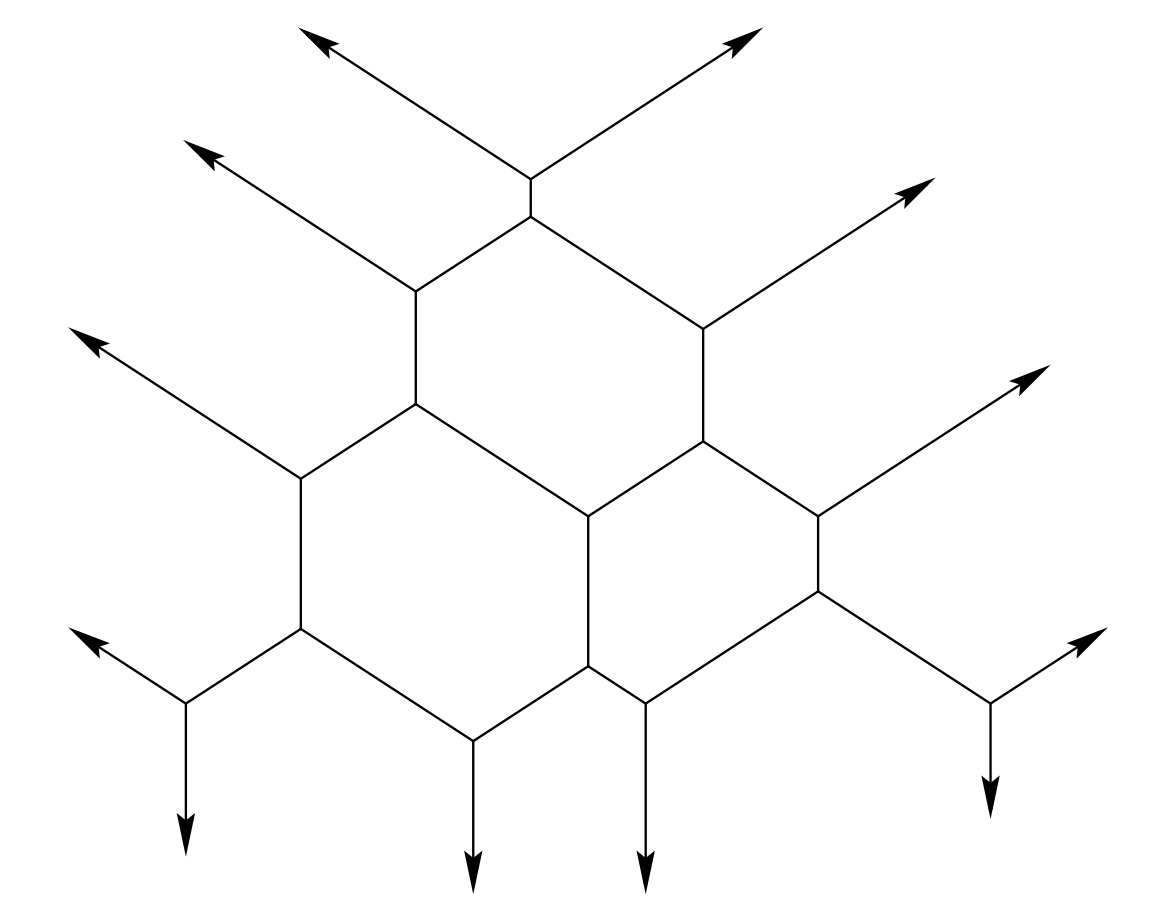}
    \caption{A honeycomb~\cite{honeycombAMS}*{Figure~1}.}
    \label{Fig:hoenycomb}
\end{figure}

An example arising somewhat further afield of geometry is the work of Allen Knutson, Terence Tao, and Christopher Woodward in representation theory~\cites{honeycombAMS, honeycomb2}.
 Knutson and Tao introduced the notion of \emph{honeycombs} (subsets of the plane as in Figure~\ref{Fig:hoenycomb}) to solve a longstanding open problem: Alfred Horn's conjectured shape of the polyhedral cone (sometimes called the Littlewood-Richardson cone) of triples of eigenvalue spectra $(\lambda, \mu, \nu)$ for Hermitian matrices $A, B, C$ which satisfy $A + B+ C = 0$~\cite{horn62}.  

\begin{figure}
    \centering
    \includegraphics[width = 1.\linewidth]{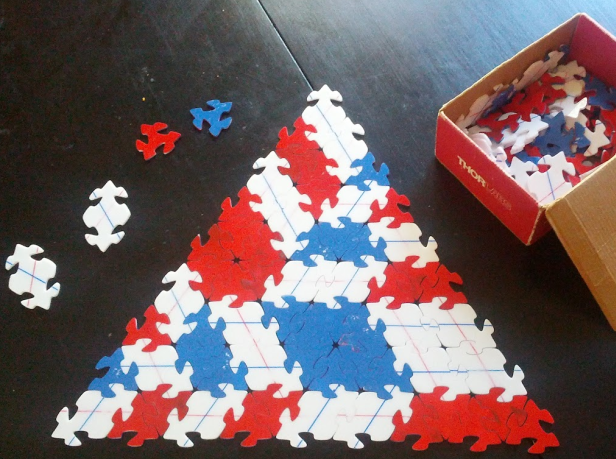}
    \caption{Knutson manufactured puzzles to study the sum-of-Hermitian-matrices problem. Image credit: Knutson.}
    \label{Fig:puzzle}
\end{figure}
 
This \emph{sum-of-Hermitian-matrices} problem has applications to perturbation theory, quantum measurement theory, and the spectral theory of self-adjoint operators. Knutson and Tao were able to show that there exist such Hermitian matrices with the specified spectra if and only if there exist honeycombs with a specified boundary. They used this correspondence to prove Horn's conjecture.  The honeycomb formalism also led naturally to a polynomial time algorithm to decide whether a triple of spectra can be realized by Hermitian matrices.  In a follow-up, Knutson, Tao, and Woodward extended the study of honeycombs to define \emph{puzzles} (Figure~\ref{Fig:puzzle}), which they described as replacing the Schubert calculus in past approaches to the Hermitian matrices problem, and used geometric arguments to give a complete characterization of the facets of the cone~\cite{honeycomb2}.  Puzzles and honeycombs provide an example of the power of rephrasing an algebraic problem as one about visual objects, where we can draw on other types of intuition.  In what circumstances can we expect these sort of insightful geometric versions to exist for algebraic problems?  When a geometric analog exists, it naturally exhibits additional features -- can we then find new corresponding objects in the original problem?  For example, what do the vertices of a honeycomb actually represent?

\begin{figure}
    \centering
    \includegraphics[width = 1.\linewidth]{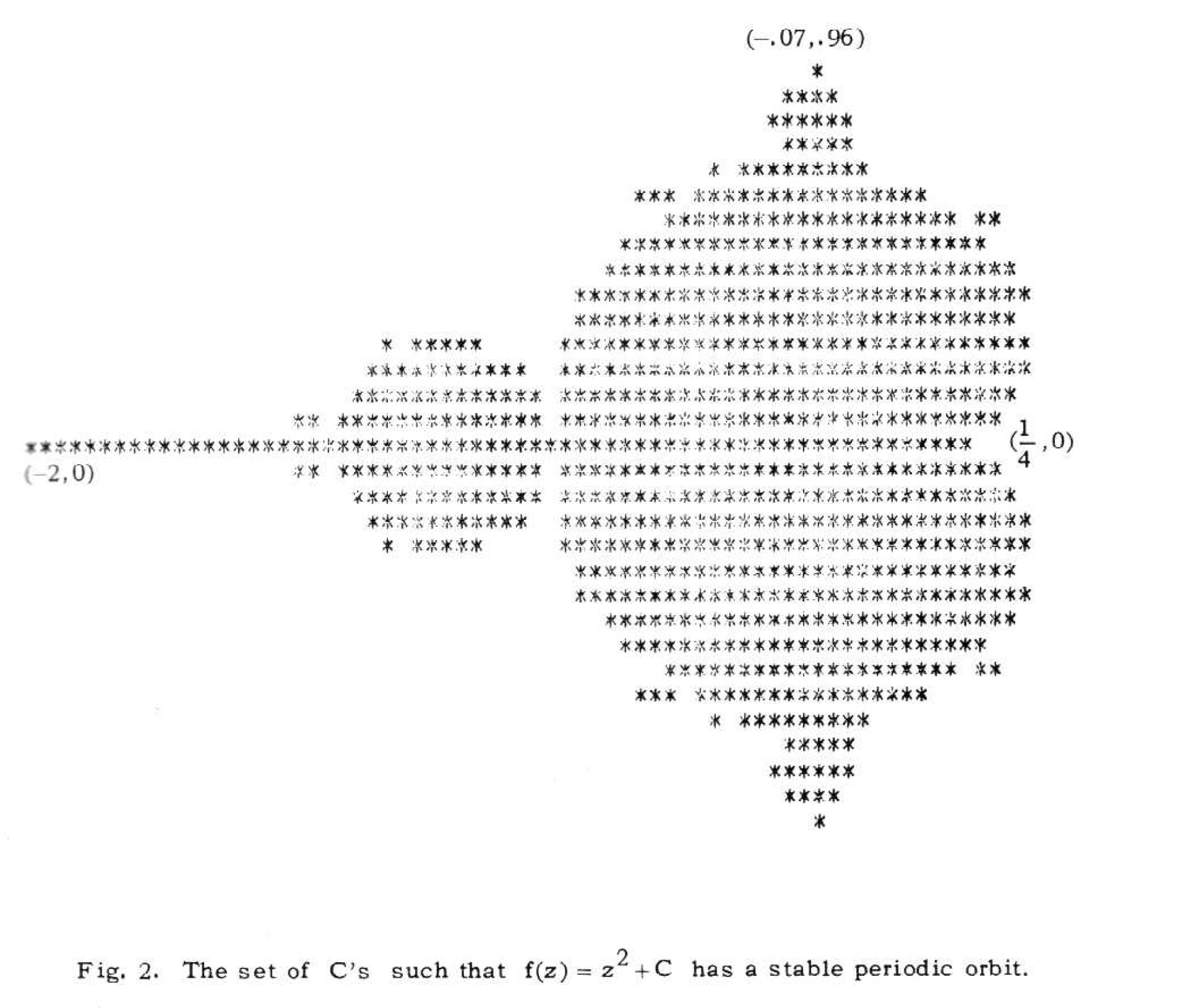}
    \caption{The first graph (in ASCII) of the Mandelbrot set. Used with permission of Princeton University Press, from Brooks and Matelski, ``The Dynamics of 2-Generator Subgroups of PSL(2, $\mathbb{C}$),'' 1980; permission conveyed through Copyright Clearance Center, Inc.~\cite{BM81}*{Figure 2}.
    }
    \label{Fig:Mandelbrot}
\end{figure}

There are, of course, many more examples.  Among these, the most famous may be the computer exploration of the Mandelbrot set and fractal geometry in the 1980's~\cite{mand82} (Figure~\ref{Fig:Mandelbrot}).  In the 1990's, Jeffrey Weeks created SnapPea (which now exists as SnapPy under the guidance of Marc Culler and Nathan Dunfield\footnote{\url{http://snappy.computop.org}}) as part of his doctoral thesis,~\cite{Weeks93}, to explore the cusp structures of hyperbolic 3-manifolds. Its use inspired David Gabai, Robert Meyerhoff, and Peter Milley to invent \emph{mom structures} to answer questions of the volumes of hyperbolic 3-manifolds~\cite{momTech}.  In the same decade, the \emph{Geometry Center} founded by Al Marden was focused on the use of computer visualization in mathematics.\footnote{\url{http://www.geom.uiuc.edu/}} It hosted mathematicians such as Eugenio Calabi, John Horton Conway, Donald E.\@ Knuth, Mumford, and Thurston, among others, and produced the GeomView software~\cite{GeomView}
used to create some famous early computer visualizations, including the sphere eversion\footnote{\textit{Outside In}, (1994), \url{http://www.geom.uiuc.edu/docs/outreach/oi/}} and illustrations for knot theory.\footnote{\textit{Not Knot}, (1991), \url{http://www.geom.uiuc.edu/video/NotKnot/}}
Illustration has shown its importance in virtually all areas of mathematics, from random tilings in combinatorics,~\cite{aztec}, 
to diagrammatic approaches to algebra,~\cite{TLDiagrams}, 
to Apollonian circle packings and Schmidt arrangements in number theory,~\cites{Mar19,Sta18,Fuchs13}, 
and their higher dimensional analogs,~\cite{spherePackings}, 
to mention just a few.

The examples above focus on pure mathematics, which is poised to join a great many other areas of scientific endeavour embracing illustration.  In applied mathematics, illustration has already made great strides.  Consider for instance the process of Alan H.\@ Schoen, when describing the \emph{gyroid} decades before it was mathematically proven to be a minimal surface.
He worked with both a sculpture of the surface and various models in \emph{Computer-Aided Design / Modelling} (CAD/CAM), which ultimately led to the structure being found in various lipid and liquid crystalline systems~\cite{schoen2012reflections}.
Other fields, like \emph{mathematical geometry processing} rely equally on quantitative measures and qualitative visualizations for judging the quality of their results~\cite{lavoue2010comparison}.
Still, a back-and-forth between the development of mathematical procedures and their application to real-world data yields results that are well-grounded in mathematical quality guarantees, yet efficient and relevant for their applications; consider for instance~\cites{alexa2022super,sawhney2022gridfree,sharp2022spelunking}.
In the field of \emph{exploratory data analysis}, visualizations even form the main tool for finding research results.
Here, large, possibly high-dimensional, data sets are investigated for patterns by embedding them, e.g., as 2D scatter plots that can then be inspected by domain experts.
With this technique, in 2020, a novel type of anti-tumor cell was discovered~\cite{vries2020highdimensional}.
None of these research results would have been possible without the utilization of illustrations.
Furthermore, this last example utilized non-linear dimensionality reduction techniques for the visualization of high-dimensional data. These techniques were themselves the result of research driven by the desire for better illustrations.

The very closely allied field of \emph{computation} in mathematics is a little ahead of illustration in its maturity as a tool for mathematical research. To give just one significant example in number theory, much recent activity has centered around the multi-million-dollar \emph{Simons Collaboration on Arithmetic Geometry, Number Theory, and Computation},\footnote{\url{https://simonscollab.icerm.brown.edu/}} whose mission states: ``Our common perspective is that advances in computational techniques accelerate research in arithmetic geometry and number theory, both as a source of data and examples, and as an impetus for effective results. The dynamic interplay between experiment, theory, and computation has historically played a pivotal role in the development of number theory.'' The work supported by the collaboration is rapidly expanding the \textit{L-Functions and Modular Forms Database},\footnote{\url{http://www.lmfdb.org}} an online database of mathematical objects (including visualizations) that is at the center of much modern progress in number theory.\footnote{See the extensive list of publications arising from the collaboration: \url{https://simonscollab.icerm.brown.edu/publications/}.}  The discipline of mathematical computation is supported by a number of journals\footnote{Consider for instance ``Advances in Computational Mathematics'', \url{https://www.springer.com/journal/10444}, the ``Journal for Computational and Applied Mathematics'', \url{https://www.sciencedirect.com/journal/journal-of-computational-and-applied-mathematics}, or the ``Journal of Computational Mathematics'', \url{https://www.jstor.org/journal/jcompmath}.} and has engendered areas of research in their own right, such as \emph{computational geometry}. 
Illustration appears to be following a similar trajectory.  As it becomes more accessible and pervasive it demands rigorous and careful study, leading to the development of mathematical illustration as a discipline in its own right.

\section*{Illustration as a discipline}

Thurston once said, ``mathematicians usually have fewer and poorer figures in their papers and books than in their heads''~\cite{thur94}.  
Although the power of good illustrations to advance mathematical knowledge is clear, they are not simple to produce.

The challenges to creating powerful and trustworthy illustrations come on many levels. 
On the one hand, some challenges are technical and concern rather practical questions regarding the production of mathematical illustrations. 
Especially with newer technologies like virtual reality or 3D modeling, the learning curves are steep and while there are general tutorials available, just a handful target issues specific to the illustration of mathematics.\footnote{A noteworthy example for introductory material, aimed at illustration of mathematics, is the \emph{Processing} tutorial of Roger Antonsen, to be found online: \url{https://rant.codes/pcmi/}.}
Consider for instance~\cite{Taalman} for a nice discussion of some of the challenges of 3D printing for mathematical illustration.

On the other hand, there are challenges within the mathematics itself.
The objects to be illustrated do not necessarily come with a description that lends itself to a suitable illustration.
Thus, a necessary initial step is the translation of the underlying mathematical object into a form that allows illustration in the first place.
However, this transformation is usually not enough by itself.
Subsequent steps aim at making the illustration effective, which can entail bridging the gap between the theoretical and the computational, 
crafting a responsive and immersive experience, or ensuring the illustration actually imparts the desired aspects of the mathematical object.
In particular the last part implies important theoretical considerations: What exactly do we want to illustrate? 
And how do we do so faithfully, i.e., without creating wrong impressions of the mathematical object illustrated?

Mathematics is not the first field of research to tackle these difficulties.  
There are parallels to be found in the development of the scientific method and statistical methods for the natural sciences: Which experimental designs and statistics can be relied upon for developing conjectures and conclusions?
Cornerstones of the scientific methods were laid down, such as the important notion of falsifiability of a scientific theory.
Similarly, statistical methods amplified their usefulness and trustworthiness when expanded from pure descriptive statistics to inference statistics and statistical tests to assert the validity of results. 
So in fact, all scientific fields have progressed by examining head-on some of the questions raised by their methodologies.

The question of \emph{illustrating well} has been asked in statistics and data visualization, as explored in Darrell Huff's best-selling book \emph{How to Lie with Statistics}, which became a standard college text.
The pioneering and richly illustrated books of Edward Tufte and Tamara Munzner on data visualization established that field in its own right. Every year, new research in data visualization is discussed at various venues, such as the Institute for Electrical and Electronics Engineers (IEEE) VIS meeting or the EuroVis conference, and published in outlets like the IEEE Transactions on Visualization and Computer Graphics. As it matures, the data visualization community addresses meta-questions on its research, such as where ``the value of visualization'' lies~\cite{van2005value} 
or ``Are we making progress in visualization research?''~\cite{correll2022we}.

Thus, the example of data visualization provides a pattern of development that  the field of mathematical illustration might follow.
However, in comparison, mathematical illustration is just taking the first steps on its journey towards being a research field. 
It is still facing basic challenges with regard to creating and evaluating the illustrations it produces.

\begin{figure*}
    \includegraphics[width=\textwidth]{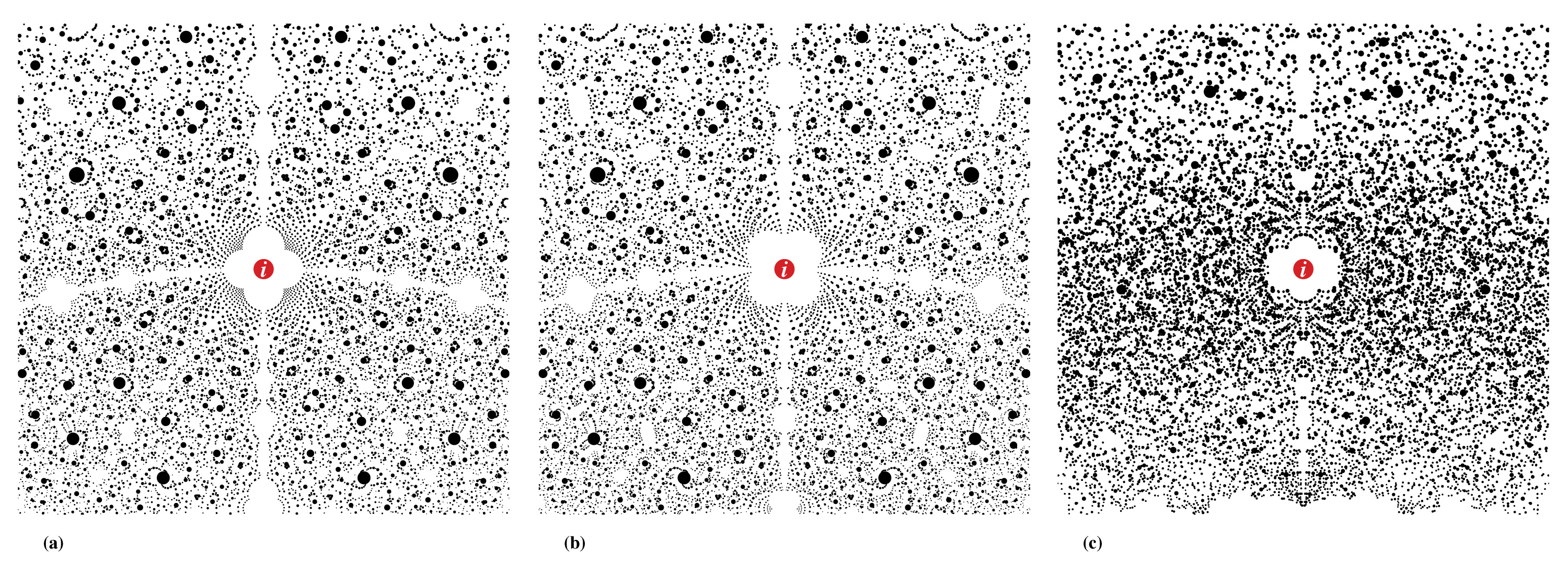}
    \caption{Roots of polynomials around $i$. Each dot is sized by $\frac{1}{\sqrt[d]{\Delta}}$, where $\Delta$ is the discriminant of the polynomial and $d$ is the degree. All images show the region of the complex plane around $i$, with real values between~$-.3$ and $.3$ and imaginary values between $.7$ and $1.3$. (a) shows all cubic roots with coefficients between $-10$ and~$10$. (b)~shows all cubic roots where the sum of the absolute value of the coefficients is less than or equal to $26$. (c)~shows all roots of polynomials up to degree $13$, with coefficients $\pm 1$~\cite{HST22} Images by Harriss, Stange, and Trettel.}
    \label{fig:aroundi}
\end{figure*}

As an example of these challenges, consider the images in Figure~\ref{fig:aroundi} showing polynomial roots near $i$ in the complex plane. 
Figure~\ref{fig:aroundi}(a) is an image of all roots of polynomials of degree $3$ with integer coefficients between $N$ and $-N$, where here $N=10$~\cite{HST22}. 
Figure~\ref{fig:aroundi}(c) is an image of all roots of polynomials with coefficients from $\{-1,1\}$ and degree no more than $D$, where in this case $D = 13$.  
In both, in the region around $i$, there appears to be a hole shaped like two ellipses overlapping at right angles. 

How to interpret this shape?  
It turns out that in Figure~\ref{fig:aroundi}(a) it is very much an artifact of the algorithm for creating these images. If you consider the picture as an approximation of all cubic roots (by allowing $N$ to tend to infinity), there are infinitely many such polynomials. By limiting $N$, we are looping through them in a growing hypercube in the coefficient space.  
The corners of this cube are the corners jutting in toward~$i$, and as the cube expands in the coefficient space, this hole will get filled in.  
If instead of looping through coefficient space in a growing cube, we choose a different ordering, the apparent shape of the void changes. 
This phenomenon is shown in Figure~\ref{fig:aroundi}(b).  

In Figure~\ref{fig:aroundi}(c), however, the limiting shape turns out to be refined but not substantially changed or filled in as $D$ increases from $13$ to infinity.
So this hole `really exists' in the picture!  
The shapes one sees at the boundaries of the limiting set of roots are explained in terms of fractal geometry and certain symmetries of this set.\footnote{These features are beautifully described by John Baez on his personal website: \url{ https://math.ucr.edu/home/baez/roots/}.} 

\begin{figure*}
    \includegraphics[width=\textwidth]{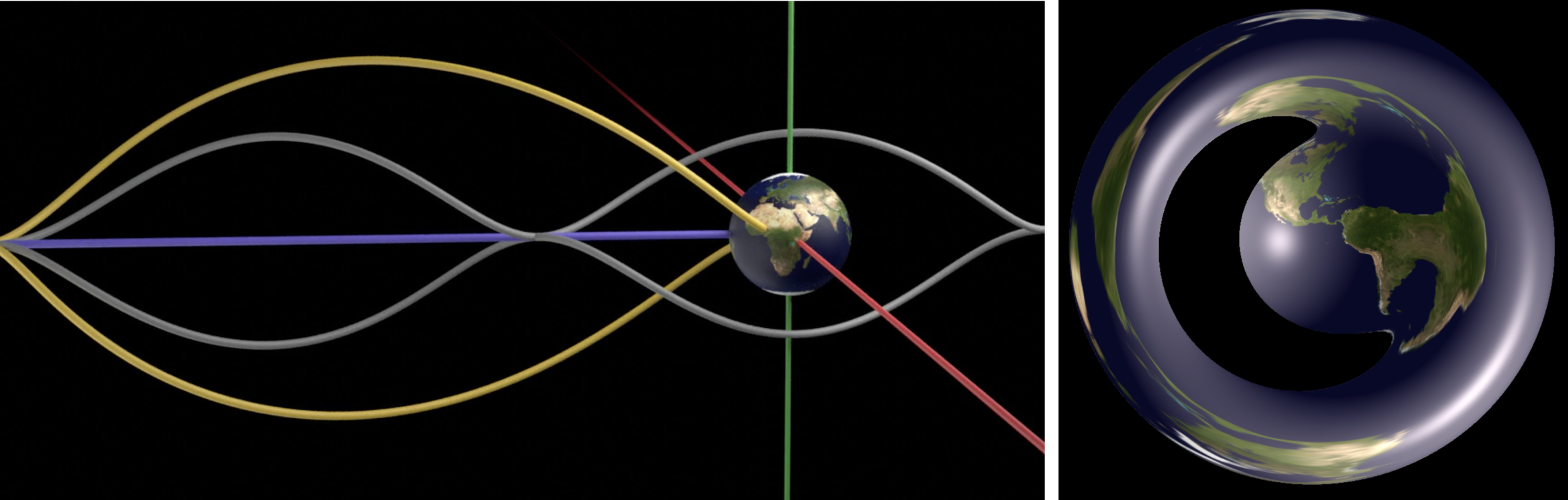}
    \caption{
    On the left: extrinsic view of some geodesics in Nil identified as a set with  $\mathbb R^3$.
    On the right: in-space view of a sphere in Nil seen from a long distance, where light is assumed to travel along geodesics.
    Images by Coulon, Matsumoto, Segerman and Trettel.
    }
    \label{fig:nil-geo}
\end{figure*}

As another example of the challenges discussed above, the virtual reality versions of Thurston's geometries of~\cite{Coulon:2020ac} are a profound way to experience these spaces, but can feel overwhelming and nearly psychedelic, as our brains struggle to make sense of what we are seeing.
As an alternative, for several of the geometries, it is possible to place the geodesics of the geometry into familiar euclidean space as curves (see Figure~\ref{fig:nil-geo}).
The interplay between these two methods of illustration can be much more enlightening than either one alone.  
The mathematical arguments that are developed to explain how one view can predict the other can end up as the basis of a mathematical proof.  
Conjectures and mathematical arguments about the space can quickly be evaluated by predicting their effect on these illustrations.

\section*{Looking forward}

Illustrations have been used both historically and in recent state-of-the-art research projects to expand the boundaries of knowledge in pure mathematics. 
Other fields of research, such as statistics and microbiology, have systematized visualization, and studied it in its own right.

However, as our gallery of examples shows, the quality of illustrations in pure mathematics varies, and there is no common framework to create, discuss, or evaluate them.
To further the possibilities that illustrations provide, there needs to be a dedicated community to tackle the next important problems.
These include, among others:
\begin{enumerate}
    \item How to identify illustrations that have rich potential to provide insight?
    \item How to identify (and mitigate) the ways that illustrations can mislead and distract?
    \item How to measure the fidelity of an illustration; are perceived patterns a result of its construction or the underlying mathematics?
    \item How can we harness the processing power and pattern-recognition capabilities of the human visual system?
    \item How can we empower a next generation of mathematical illustrators to create and leverage sophisticated illustrations?
    \item And how do we increase professional recognition of the illustration of mathematics?
\end{enumerate}
Exploring these questions will lay the foundation of a discipline built around the illustration of mathematics, providing powerful tools for the advancement of mathematical research.

\subsubsection*{Acknowledgment}

The authors would like to thank the two referees for their valuable feedback and suggestions. The illustrating mathematics group also provided insightful discussions on the topics presented here.

\bibliographystyle{plain}
\bibliography{bibliography}

\end{document}